\documentclass[11pt]{jdg-p}

\usepackage{amsmath,amscd}
\usepackage{amsthm}
\usepackage{amssymb}
\usepackage{amsfonts}
\usepackage{mathrsfs}
\usepackage{latexsym}
\usepackage{graphicx}
\usepackage{epsfig}
\usepackage[T1]{fontenc}

\theoremstyle{plain}
\newtheorem{theorem}{Theorem}[section]
\newtheorem{lemma}[theorem]{Lemma}

\newtheorem{corollary}[theorem]{Corollary}

\newtheorem*{claim}{Claim}

\theoremstyle{definition}
\newtheorem{definition}[theorem]{Definition}

\theoremstyle{remark}
\newtheorem{remark}{Remark}
\newtheorem*{acknowledgement}{Acknowledgement} 

\renewcommand{\proofname}{\bf{Proof}}

\renewcommand{\qed}{\mbox{}\hfill\openbox}

\renewcommand{\labelenumi}{{\textrm (\roman{enumi})}}

\numberwithin{equation}{section}

\allowdisplaybreaks

\def\XXint#1#2#3{{\setbox0=\hbox{$#1{#2#3}{\int}$}
    \vcenter{\hbox{$#2#3$}}\kern-.5\wd0}}

\renewcommand{\leq}{\leqslant}
\renewcommand{\geq}{\geqslant}

\newcommand{\inner}[2]{\langle #1\,,#2\rangle}

\newcommand{\C}{\mathbb{C}}
\newcommand{\R}{\mathbb{R}}
\renewcommand{\H}{\mathbb{H}}
\renewcommand{\L}{\mathbb{L}}

\newcommand{\T}{\mathbb{T}}

\newcommand{\Acal}{\mathcal{A}}
\newcommand{\Bcal}{\mathcal{B}}
\newcommand{\Hcal}{\mathcal{H}}
\newcommand{\Ocal}{\mathcal{O}}

\newcommand{\Cscr}{\mathscr{C}}
\newcommand{\Lscr}{\mathscr{L}}
\newcommand{\Nscr}{\mathscr{N}}

\DeclareMathOperator{\PSL}{PSL}

\DeclareMathOperator{\Area}{Area}
\DeclareMathOperator{\Isom}{Isom}
\DeclareMathOperator{\Length}{Length}
\DeclareMathOperator{\Vol}{Vol}

\DeclareMathOperator{\diam}{diam}
\DeclareMathOperator{\dist}{dist}

\renewcommand{\epsilon}{\varepsilon}

\newcommand{\gbar}{\bar{g}}

\begin{document}

\title{Minimal Surfaces in Quasi-Fuchsian 3-Manifolds}

\author{Biao Wang}
\date{March 28, 2009}
\address{Department of Mathematics\\
         Cornell University\\
         Ithaca, NY 14853-4201\\}
\email{wang@math.cornell.edu}

\begin{abstract}
In this paper, we prove that if a quasi-Fuchsian $3$-manifold 
$M$ contains a closed geodesic with complex length 
$\Lscr=l+i\theta$ such that $\theta/l\gg{}1$, then it contains 
at least two minimal surfaces which are incompressible in $M$.
\end{abstract}

\maketitle


\section{Introduction}

For any complete hyperbolic $3$-manifold $N$, there is a 
holonomy representation $\rho:\pi_{1}(N)\to\Isom^{+}(\H^{3})$
so that $G:=\rho(\pi_{1}(N))$ is a discrete subgroup of
$\Isom^{+}(\H^{3})$ called {\em Kleinian group} and 
$N=\H^{3}/G$. For any torsion free Kleinian group 
$G\subset\Isom^{+}(\H^{3})$, the quotient 
$\H^{3}/G$ is a complete hyperbolic $3$-manifold. 

For a torsion free Kleinian group $G$, if its limits set 
$\Lambda(G)\subset{}S_{\infty}^{2}$ is a closed Jordan 
curve, i.e. $\Lambda(G)$ is homeomorphic to the unit circle
$S^{1}$, then $G$ is called a {\em quasi-Fuchsian group}, and
$M=\H^{3}/G$ is called a {\em quasi-Fuchsian $3$-manifold}.
The quasi-Fuchsian $3$-manifold $M$ has a {\em convex core},
denoted by $\Cscr(M)$, which is the quotient of the convex 
hull of $\Lambda(G)$.
Topologically, $M$ is homeomorphic to $\Sigma\times\R$, where
$\Sigma$ is a finite type Riemann surface with negative
Euler characteristic. In this paper, $\Sigma$ is assumed to
a closed surface with genus $\geq{}2$ except in 
section \ref{sec:coarea} where $\Sigma$ 
denotes a hypersurface in a Riemannian manifold. 

Any quasi-Fuchsian $3$-manifold $M=\Sigma\times\R$ 
always contains a minimal surface, say $T$, which is
incompressible in $M$ 
(cf. \cite{anderson1983,SY1979(annals),uhlenbeck1983}). 
Besides, K. Uhlenbeck 
also proved that if the principle curvature of $T$ is between 
$-1$ and $1$, then $T\subset{}M$ is unique
and $M$ can be foliated by equidistant surfaces to $T$
(cf. \cite{uhlenbeck1983}). 
In this case, we proved that $M$ admits a foliation such 
that each leaf is a surface of constant mean curvature 
(cf. \cite{wang2008}).

In this paper, we will prove the following theorem.

\begin{theorem}\label{thm:main theorem}
Let $M=\Sigma\times\R$ be a quasi-Fuchsian $3$-manifold, 
let $\gamma\subset{}M$ be a simple closed geodesics $\gamma$ 
with complex length $\Lscr=l+i\theta$, and let
$S\subset{}M\setminus\gamma$ be an embedded closed surface
that is incompressible in $M$. If $\theta/l\gg{}1$,
then there exists a least area minimal surface
$T\subset{}M\setminus\gamma$ that is isotopic
to $S$ in $M\setminus\gamma$, here $T$ is of least area means 
if there is another minimal surface 
$T'\subset{}M\setminus\gamma$ isotopic
to $S$, then $\Area(T)\leq\Area(T')$.
\end{theorem}

If a quasi-Fuchsian $3$-manifold  
contains more than one sufficiently short simple closed 
geodesics, say $\gamma_{1},\ldots,\gamma_{n}$, then we have 
a result similar to Theorem \ref{thm:main theorem}, see
Corollary \ref{cor:corollary of main theorem} for detail. 

The plan of this paper is as follows.
In section \ref{sec:tube} we review the basic properties
of tubes, in section \ref{sec:helicoid} we discuss some
properties of helicoids in $\H^{3}$, in section \ref{sec:coarea} 
we prove one version of coarea formulae used frequently 
in \cite{CG2006}, and in section \ref{sec:shrinkwrapping} we 
introduce a tool called shrinkwrapping developed by
Calegari and Gabai in \cite{CG2006}, in section 
\ref{sec:proof} we prove Theorem \ref{thm:main theorem}.

\begin{acknowledgement} 
I would like to express my
gratitude to my advisor Profesoor Bill Thurston for his
help, stimulating suggestions and encouragement in the past
several years.
\end{acknowledgement}

\section{Tubes of short geodesics}\label{sec:tube}

We consider the upper half space model of hyperbolic $3$-space,
that is, a three-diemnsional space
\begin{equation*}
   \H^{3}=\{z+tj\ |\ z\in\C,\ t\in\R,\ t>0\}
\end{equation*}
which is equipped with the (hyperbolic) metric
\begin{equation*}
   ds^{2}=\frac{|dz|^{2}+dt^{2}}{t^{2}}\ ,
\end{equation*}
where $z=x+iy$ for $x,y\in\R$. Let $\PSL_{2}(\C)$ denote the 
set of M{\"o}bius transformations on $\H^{3}$. For any 
$g\in\PSL_{2}(\C)$, write
\begin{equation*}
   g(z)=\frac{az+b}{cz+d}\ ,
   \quad\forall\,z\in\C\ ,
\end{equation*}
then the Poincar{\'e} extension of $g$ is given by
\begin{equation}
   g(z+tj)=\frac{(az+b)(\overline{cz+d})+a\overline{c}t^{2}+tj}
   {|cz+d|^{2}+|c|^{2}t^{2}}\ .
\end{equation}
In this paper, we are interested in the loxodromic elements,
i.e. by conjugacy, $g(z)=\alpha{}z$, where
$\alpha=\exp(l+i\theta)$ with $l>0$ and $0\leq\theta<2\pi$.
The loxodromic transformation
\begin{equation}
   g(z+tj)=\alpha{}z+|\alpha|tj
\end{equation}
translates points on the $t$-axis by distance $l$ and
twists a normal plane by angle $\theta$. The complex number
$\Lscr:=l+i\theta$ is called the {\em complex length} of 
$g$. Since the trace of $g$ is invariant under conjugate
transformations, the complex length of $g$ is also invariant
under conjugate transformations. Notice that $g$ determines 
$\theta$ only up to a multiple of $2\pi$.

Now suppose $\gamma$ is a simple closed geodesic in a 
complete hyperbolic $3$-manifold $N$, then a lift 
$\tilde\gamma$ in $\H^{3}$ is the axis of a loxodromic element
$g\in\PSL_{2}(\C)$ representing $\gamma$.
The choice of different 
lift $\tilde\gamma{}'$ gives rise to another element of
$\PSL_{2}(\C)$ such that two elements are conjugate in
$\PSL_{2}(\C)$. Two lifts $\tilde\gamma$ and
$\tilde\gamma{}'$  possesses a common perpendicular.
Let $r_{\max}(\gamma)$ be the maximal number $r>0$ such that
$\Nscr_{r}(\tilde\gamma{})\cap\Nscr_{r}(\tilde\gamma{}')=\emptyset$
as $\tilde\gamma{}'$ runs through all the lifts of $\gamma$
different from $\tilde\gamma$, where 
$\Nscr_{r}(\tilde\gamma)=\{x\in\H^{3}\ |\
\dist(x,\tilde\gamma)<r\}$ 
is the $r$-neighborhood of the geodesic $\tilde\gamma$ 
in $\H^{3}$.  We call $r_{\max}=r_{\max}(\gamma)$ the 
{\em tube radius} of $\gamma$. The (maximal solid) tube of 
$\gamma$ is defined by
\begin{equation}\label{eq:solid tube}
   \T(\gamma)=\Nscr_{r_{\max}}(\tilde\gamma)/
   \langle{}g\rangle\ .
\end{equation}

When $l$, the real length of $\gamma$, is sufficiently small,
R. Meyerhoff found a relation between $r_{\max}$ and $l$.

\begin{theorem}[\cite{meyerhoff1987}]
\label{thm:meyerhoff1987-(4)}
If $\gamma$ is a closed geodesic in a complete hyperbolic
$3$-manifold $N$, with complex length $\Lscr=l+i\theta$
satisfying
\begin{equation}\label{eq:upper bound of length}
   l\leq\frac{\sqrt{3}}{4\pi}\,[\log(\sqrt{2}+1)]^{2}
   \approx{}0.107\ ,
\end{equation}
then there exists an embedded solid tube around $\gamma$ whose
radius is given by the following formula
\begin{equation}
   \sinh^{2}r_{\max}(\gamma)=\frac{1}{2}\left(
   \frac{\sqrt{1-2\kappa(l)}}{\kappa(l)}-1\right)\ ,
\end{equation}
where $\kappa(l)=\cosh(\sqrt{4\pi{}l}/\sqrt[4]{3})-1$.
\end{theorem}

If a hyperbolic $3$-manifold $N$ contains more than one 
short geodesics, we want the solid tubes of different 
geodesics to be disjoint. Fortunately, R. Meyerhoff also
proved the following theorem.

\begin{theorem}[\cite{meyerhoff1987}]
\label{thm:meyerhoff1987-(7)}
The solid tubes constructed by \eqref{eq:solid tube}
about different short geodesics whose real length
satisfy \eqref{eq:upper bound of length} do not intersect.
\end{theorem}

\section{Helicoids in $\H^{3}$}\label{sec:helicoid}

In this section, we consider the Lorentzian $4$-space
$\L^{4}$, i.e. a vector space $\R^{4}$ with the Lorentzian
inner product
\begin{equation*}
   \inner{x}{y}=-x_{1}y_{1}+x_{2}y_{2}+x_{3}y_{3}+x_{4}y_{4}
\end{equation*}
where $x,y\in\R^{4}$. The hyperbolic space $\H^{3}$ can be
considered as the unit sphere of $\L^{4}$:
\begin{equation*}
   \H^{3}=\{x\in\L^{4}\ |\ \inner{x}{x}=-1,\,x_{1}\geq{}1\}\ .
\end{equation*}
The helicoid $\Hcal_{a}$ is the surface parametrized by the
$(u,v)$-plane in the following way: 
\begin{equation}
   \Hcal_{a}=\left\{x\in\H^{3}\ \left|
   \begin{aligned}
      &x_{1}=\cosh{}u\cosh{}v, &&x_{2}=\cosh{}u\sinh{}v\\
      &x_{3}=\sinh{}u\cos(av), &&x_{4}=\sinh{}u\sin(av)
   \end{aligned}
   \right.\right\}
\end{equation}
The axis of the helicoid $\Hcal_{a}$ is
\begin{equation*}
   (\cosh{}v,\sinh{}v,0,0)\ ,
   \quad{}-\infty<v<\infty\ ,
\end{equation*}
which is the intersection of $x_{1}x_{2}$-plane and the
three dimensional hyperboloid. The first fundamental
form of $\Hcal_{a}$ is
\begin{equation}
   I=du^{2}+(\cosh^{2}u+a^{2}\sinh^{2}u)dv^{2}\ .
\end{equation}

Intuitively, when $a>0$ is sufficiently large, the helicoid
$\Hcal_{a}$ will become unstable. Actually, Hiroshi Mori
proved the following theorem.

\begin{theorem}[\cite{mori1982}]
Let $\Hcal_{a}$ be the helicoid as above, then it
is a minimal surface in $\H^{3}$. Furthermore,
\begin{itemize}
   \item $\Hcal_{a}$ is globally stable if
         $a\leq{}3\sqrt{2}/4$,
   \item $\Hcal_{a}$ is globally unstable if
         $a\geq\sqrt{105\pi}/8$.
\end{itemize}
\end{theorem}

The unstability of $\Hcal_{a}$ when $a$ is sufficiently large
is very important in this paper. We also need the following
lemma, which shows that $\Hcal_{a}$ is a ruled surface in
$\H^{3}$.

\begin{lemma}[Tuzhilin\ \cite{tuzhilin1993}]
\label{lem:tuzhilin1993}
Let $\Hcal_{a}$ be the helicoid as above, then it is a ruled
surface, i.e. it is stratified into straight lines in the sense
of hyperbolic metric.
\end{lemma}

Let $\gamma$ be a closed geodesic in a hyperbolic
$3$-manifold $N$, let $g\in\PSL_{2}(\C)$ be a loxodromic
element with complex length $\Lscr=l+i\theta$ which represents
$\gamma$, and let $\tilde\gamma\subset\H^{3}$ be the axis of $g$.
Without of lose generality, we may assume $\theta>0$. Write
$a=\theta/l$, let $\Hcal_{\theta/l}\subset\H^{3}$ be the
helicoid constructed as above, and let
\begin{equation}
   \Acal_{\theta/l}=\frac{\Hcal_{\theta/l}
   \cap\Nscr_{r_{\max}}(\tilde\gamma)}
   {\langle{}g\rangle}\ .
\end{equation}

\begin{lemma}\label{lem:area-of-annulus}
The quotient surface $\Acal_{\theta/l}$
is a minimal annulus in the tube $\T(\gamma)$, and its area is
\begin{equation}
   \Area(\Acal_{\theta/l})=
   2\int_{0}^{r_{\max}}(l^{2}\cosh^{2}u+
   \theta^{2}\sinh^{2}u)^{1/2}du\ .
\end{equation}
\end{lemma}

\begin{proof}For each point $p\in\Hcal_{\theta/l}$, the point
$g(p)\in\Hcal_{\theta/l}$ is obtained by twisting by angle
$\theta$ and translating by distance $l$, so
$\Hcal_{\theta/l}$ is invariant under the cyclic group
$\langle{}g\rangle$, the quotient surface $\Acal_{\theta/l}$
is an annulus in $\T(\gamma)$. By Lemma \ref{lem:tuzhilin1993}
and coarea formula \eqref{eq:coarea formula II}, the area of
$\Acal_{\theta/l}$ can be computed as follows:
\begin{align*}
   \Area(\Acal_{\theta/l})
         &=\int_{0}^{r_{\max}}\Length(\Acal_{\theta/l}
           \cap\Nscr_{r}(\gamma))dr\\
         &=2\int_{0}^{r_{\max}}\left\{
            \int_{0}^{l}\left[\cosh^{2}r+
            \left(\frac{\theta}{l}\right)^{2}
            \sinh^{2}r\right]^{1/2}ds\right\}dr\\
         &=2\int_{0}^{r_{\max}}(l^{2}\cosh^{2}r+
            \theta^{2}\sinh^{2}r)^{1/2}dr\ .
\end{align*}
Therefore the lemma is proved.
\end{proof}

\begin{theorem}\label{thm:area comparison}
Suppose a closed geodesic $\gamma$ in a complete
hyperbolic $3$-manifold $N$ has complex length $\Lscr=l+i\theta$,
let $\T(\gamma)$ and $\Acal_{\theta/l}$ be constructed as above.
If $\theta/l$ is sufficiently large, then
\begin{equation*}
   \Area(\Acal_{\theta/l})>\Area(\partial\T(\gamma))\ .
\end{equation*}
\end{theorem}

\begin{proof}By Lemma \ref{lem:area-of-annulus}, when
$\theta/l\gg{}1$, we have the following estimate
\begin{align*}
   \Area(\Acal_{\theta/l})
      \geq{}2\theta\int_{0}^{r_{\max}}\sinh{u}\,du
      =2\theta(\cosh{}r_{\max}-1)\ .
\end{align*}
Since $r_{\max}(\gamma)\to\infty$ as $l\to{}0$, we have
the estimate
\begin{equation*}
   \Area(\Acal_{\theta/l})\geq\theta\cosh{}r_{\max}
\end{equation*}
if $l\ll{}1$. On the other hand, the area of the torus 
$\partial\T(\gamma)$ is given by
\begin{equation*}
   \Area(\partial\T(\gamma))=2\pi{}l\cosh{}r_{\max}
   \sinh{}r_{\max}\ .
\end{equation*}
As $l\to{}0$, we claim that $l\sinh{}r_{\max}\to{}0$.
In fact, by Theorem \ref{thm:meyerhoff1987-(4)}, we have 
the identity
\begin{equation*}
   l^{2}\sinh^{2}(r_{\max})=\frac{l^{2}}{2}
   \left(\frac{\sqrt{1-2\kappa(l)}}{\kappa(l)}-1\right)\ ,
\end{equation*}
where $\kappa(l)=\cosh(\sqrt{4\pi{}l}/\sqrt[4]{3})-1=
2\pi{}l/\sqrt{3}+\Ocal(l^{2})$. Therefore
\begin{equation*}
   l^{2}\sinh^{2}(r_{\max})=
   \frac{l}{2\pi/\sqrt{3}+\Ocal(l)}\cdot
   \frac{\sqrt{1-2\kappa(l)}-\kappa(l)}{2}
   \rightarrow{}0\ ,
\end{equation*}
as $l\to{}0$. Therefore when $\theta/l\gg{}1$, we 
have the estimate $2\pi{}l\sinh{}r_{\max}<\theta$,
which implies the lemma.
\end{proof}

\begin{remark}The area of the meridian disk of $\T(\gamma)$
is $2\pi(\cosh{}r_{\max}-1)$, so the computation in
Theorem \ref{thm:area comparison} also implies that the
area of the meridian disk is bigger than 
$\Area(\partial\T(\gamma))$ when $\theta/l$ is
sufficiently large.
\end{remark}

\begin{remark}From the computation above, as $l\to{}0$, we have
\begin{equation*}
   \Area(\partial\T(\gamma))\rightarrow\frac{\sqrt{3}}{2}
   \qquad\text{and}\qquad
    \Vol(\T(\gamma))\rightarrow\frac{\sqrt{3}}{4}\ .
\end{equation*}
\end{remark}


\section{Coarea formula}\label{sec:coarea}

The following coarea formulae \eqref{eq:coarea formula I}
and \eqref{eq:coarea formula II} were used frequently in
\cite{CG2006} without proof. Since the coarea formulae are very
important in this paper too, we give a proof here for convenience,
which is also pretty simple.

\begin{lemma}\label{lem:coarea formula}
Suppose $\Sigma$ is a hypersurface in a complete
Riemannian manifold $(N,\gbar_{ij})$. Let $p\in{}N$ be a fixed
point, for any point $q\in\Sigma$, define $\alpha(q)$ to be
the angle between the tangent space to $\Sigma$ at $q$, and
the radial geodesic through $q$ emanating from $p$. Then
\begin{equation}\label{eq:coarea formula I}
   \Area(\Sigma\cap\Nscr_{s}(p))=
   \int_{0}^{s}\int_{\Sigma\cap\partial\Nscr_{t}(p)}
   \frac{1}{\cos\alpha}\,dldt\ .
\end{equation}
\end{lemma}

\begin{proof}Define a function $h:\Sigma\to\R$ by
$h(q)=\dist(p,q)$, where $\dist(\cdot,\cdot)$ is the distance
function on $N$. Then, by the standard coarea formula
(cf. \cite[p. 89]{SY1994(book)}), we have
\begin{equation*}
   \Area(\Sigma\cap\Nscr_{s}(p))=
   \int_{0}^{s}\int_{\Sigma\cap\partial\Nscr_{t}(p)}
   \frac{1}{|\nabla_{\Sigma}h|}\,dldt\ .
\end{equation*}
Since $\nabla{}h=(\nabla{}h)^{\bot}+(\nabla{}h)^{\parallel}=
\nabla^{\bot}h+\nabla_{\Sigma}h$
and $|\nabla{}h|_{\gbar_{ij}}=1$, one has
$|\nabla_{\Sigma}h|=\cos\alpha$.
\end{proof}

\begin{remark}If $\gamma$ is a simple closed geodesic in $N$,
then we still have the following coarea formula
\begin{equation}\label{eq:coarea formula II}
   \Area(\Sigma\cap\Nscr_{s}(\gamma))=
   \int_{0}^{s}\int_{\Sigma\cap\partial\Nscr_{t}(\gamma)}
   \frac{1}{\cos\alpha}\,dldt\ .
\end{equation}
\end{remark}

\section{Shrinkwrapping}\label{sec:shrinkwrapping}

In this section, we use a technical tool, called
shrinkwrapping developed by Calegari and Gabai
in \cite{CG2006}, to find minimal surfaces in quasi-Fuchsian
$3$-manifolds containing a very short simple  closed geodesic.
Roughly speaking,
given a simple closed geodesic $\gamma$ in a quasi-Fuchsian
$3$-manifold $M$ and an embedded surface
$S\subset{}M\setminus\gamma$, a surface $T\subset{}M$ is
obtained from $S$ by {\em shrinkwrapping $S$ rel. $\gamma$}
if it homotopic to $S$, can be approximated by an isotopy from
$S$ supported in $M\setminus\gamma$, and is least area
subject to these constraints.

\begin{definition}\label{def:CG2006-(1.17)}
Choose some positive $\sigma$ so that 
$\Nscr_{\sigma}(\gamma)\subset{}M$
is an embedded solid torus , for $0\leq{}t<1$ we
define a family of Riemannian metrics $(M,g_{t})$ in the
following manner. 

Let $h:\Nscr_{\sigma(1-t)}(\gamma)\to[0,\sigma(1-t)]$
be the function whose value at a point $p$ is the hyperbolic
distance from $p$ to $\gamma$. Let $\phi:[0,1]\to[0,1]$ be a
$C^{\infty}$ bump function given by
   \begin{equation*}
       \phi(t)=
       \begin{cases}
          0, & t\in[0,1/4],\\
          \text{strictly increasing}, & t\in[1/4,1/3],\\
          1, & t\in[1/3,2/3],\\
          \text{strictly decreasing}, & t\in[2/3,3/4],\\
          0, & t\in[3/4,1]\ .
       \end{cases}
   \end{equation*}
We define a metric $g_{t}$ on $M$:
\begin{itemize}
   \item $g_{t}$ agrees with the hyperbolic metric on
         $M\setminus\Nscr_{\sigma(1-t)}(\gamma)$;
   \item $g_{t}$ is conformally equivalent to the hyperbolic
         metric on $\Nscr_{\sigma(1-t)}(\gamma)$, which is more
         precisely given by the ratio
         \begin{equation*}
            \frac {g_{t}\ \text{length element}}
            {\text{hyperbolic length element}} =
            1 + 2\phi\left(\frac {h(p)}{\sigma(1-t)}\right)\ .
         \end{equation*}
\end{itemize}
\end{definition}

\begin{lemma}[\cite{CG2006}]\label{lem:CG2006-(1.18)}
The $g_{t}$ metric satisfies the following properties:
   \begin{enumerate}
      \item For each $t$ there is an $f(t)$ satisfying
            $\frac{2}{3}\,\sigma(1-t)<f(t)<\frac{3}{4}\,\sigma(1-t)$
            such that the union of tori
            $\partial\Nscr_{f(t)}(\gamma)$ are totally
            geodesic for the $g_{t}$ metric.
      \item The restricted metric
            $g_{t}|{\Nscr_{\sigma}(\gamma)}$
            admits a family of isometries preserving
            $\gamma$ and acts transitively on
            the unit normal bundle (in $M$) to $\gamma$.
     \item The area of a disk cross--section on
           $\Nscr_{\sigma(1-t)}$ is $O((1-t)^2)$.
     \item The metric $g_{t}$ dominates the hyperbolic metric on
           $2$-planes. That is, for all $2$-vectors $\nu$, the
           $g_{t}$ area of $\nu$ is at least as large as the
           hyperbolic area of $\nu$.
   \end{enumerate}
\end{lemma}

\begin{lemma}[\cite{CG2006}]\label{lem:CG2006-(1.20)}
Let $M$, $\gamma$ and $S$ be as in the statement of
Theorem \ref{thm:main theorem}. Let $f(t)$ be a function
given in Lemma \ref{lem:CG2006-(1.18)} so that
$\partial\Nscr_{f(t)}(\gamma)$ is totally geodesic
with respect to the $g_{t}$ metric. Then for each $t$, there
exists an embedded surface $S_{t}$ isotopic in
$M\setminus\Nscr_{f(t)}(\gamma)$ to $S$, and which
is globally $g_{t}$-least area among all such surfaces.
\end{lemma}

\begin{proof}By Lemma \ref{lem:CG2006-(1.18)}, the surfaces
$\partial\Nscr_{f(t)}(\gamma)$ are
totally geodesic with respect to the $g_{t}$ metrics,
therefore these surfaces can act as barrier surfaces. As
a quasi-Fuchsian $3$-manifold, $M$ has a compact convex core
$\Cscr(M)$ whose boundary consists of two incompressible
surfaces which are convex with respect to the inward normal. Let
\begin{equation*}
   M_{t}=\Cscr(M)\setminus\Nscr_{f(t)}(\gamma)\ ,
\end{equation*}
then $M_{t}$ is a compact manifold with mean convex boundary.
We assume, after a small isotopy if necessary, that $S$ does
not intersect $\Nscr_{f(t)}(\gamma)$ for any $t\in[0,1]$,
therefore $S$ is an incompressible surface in $M_{t}$.
Then Meeks-Simon-Yau's main theorem in \cite{MSY1982}
implies that there exists a globally least area surface
$S_{t}$ with respect to the metric $g_{t}$.
\end{proof}

\section{Proof of Theorem \ref{thm:main theorem}}
\label{sec:proof}

\begin{lemma}\label{lem:MS-short-curve}
Suppose $\gamma$ is a closed geodesic in a
complete hyperbolic $3$-manifold $M$, let $C$ be a simple
null-homotoptic curve on $\partial\Nscr_{r}(\gamma)$ with
$\Length(C)<{}2\pi\sinh{}r$, here $r>0$ is chosen so that the
$r$-tube $\Nscr_{r}(\gamma)$ is embedded, then $C$ bounds an
embedded least area disk-type minimal surface
$\Delta\subset\Nscr_{r}(\gamma)\setminus\gamma$. 
\end{lemma}

\begin{proof}Suppose $C$ bounds a disk
$D\subset\partial\Nscr_{r}(\gamma)$, then the diameter of $D$,
denoted by $\diam(D)$, with respect to the induced metric on
$\partial\Nscr_{r}(\gamma)$ is $<\pi\sinh{}r$; otherwise we have
$\Length(C)\geq{}2\pi\sinh{}r$, which implies a contradiction.

Since the longitudes
$\Acal_{\theta/l}\cap\partial\Nscr_{r}(\gamma)$ divide the torus
$\partial\Nscr_{r}(\gamma)$ into two semi-tori $R_{\pm}$, and recall
that the length of the meridian of $\partial\Nscr_{r}(\gamma)$ is
$2\pi\sinh{}r$, so $D$ must be contained in one of the semi-tori,
say $R_{+}$. These longitudes bound a minimal annulus
$\Acal_{\theta/l}(r):=\Acal_{\theta/l}\cap\Nscr_{r}(\gamma)$, 
then $R_{+}\cup\Acal_{\theta/l}(r)$ is
mean convex with respect to the inward normal, by a result in
\cite{MY1982(t),MY1982(mz)}, we know that $C$ bounds a least area
disk $\Delta$ contained in the semi-tube bounded by $R_{+}$ and
$\Acal_{\theta/l}(r)$, and any such least area disk is properly
embedded. The minimal disk $\Delta$ cannot intersect
$\Acal_{\theta/l}(r)$, since the latter is a barrier surface. In
particular, $\Delta$ must be disjoint from the closed geodesic
$\gamma\subset\Acal_{\theta/l}(r)$.
\end{proof}

\begin{remark}If $\gamma$ is very short and 
$r\in[r_{\max}/4,r_{\max}/2]$, where $r_{\max}$ is the tube 
radius of $\gamma$, then Lemma \ref{lem:MS-short-curve}
can be proved without using Meeks-Yau's result.
\begin{itemize} 
\item The existence of $\Delta$ is from J. Douglas, T. Rad{\'o}
      and C. B. Morrey.
\item If $\Delta$ intersects $\gamma$, then $\Delta$ contains an 
      intrinsic disk $\Bcal_{r}$ of radius $r$, and 
      $\Area(\Bcal_{r})\geq{}2\pi(\cosh{}r-1)$. Since $D$ must 
      be contained in a semi-torus, we have
      $\Area(D)\leq\pi{}l\sinh{}r\cosh{}r$. Therefore 
      $2\pi(\cosh{}r-1)\leq\pi{}l\sinh{}r\cosh{}r$. This is 
      impossible when $r_{\max}$ is sufficiently large.
\item If $\Delta\not\subset\Nscr_{r}(\gamma)$, then the same 
      argument as above tells us that
      $\Delta\subset\Nscr_{2r}(\gamma)$,
      and so there exists $r_{0}\in(r,2r)$ such that 
      $\partial\Nscr_{r_{0}}(\gamma)$ touches $\Delta$ for the 
      first time from outside, but this is impossible by maximum
      principle.
\end{itemize}
Therefore $C$ bounds a least area minimal disk 
$\Delta\subset\Nscr_{r}(\gamma)\setminus\gamma$.
\end{remark}

Now we can prove Theorem \ref{thm:main theorem}.

\begin{proof}[\indent{\bf Proof of Theorem \ref{thm:main theorem}}]
Let $\sigma<r_{\max}/2$ be a constant, where $r_{\max}$ is the
tube radius of $\gamma$ that will be determined later.
By Lemma \ref{lem:CG2006-(1.20)}, we can construct
$g_{t}$-least area minimal surface $S_{t}$ that is isotopic
to $S$ for $t\in[0,1)$.

\begin{claim}Suppose $\theta/l\gg{}1$, then
there exists a constant $\epsilon=\epsilon(l,\theta)$
such that $S_{t}\cap\Nscr_{\epsilon}(\gamma)=\emptyset$ for
all $t$ that is sufficiently large.
\end{claim}

\begin{proof}[{\indent\bf Proof of Claim}]At first, by
Lemma \ref{lem:MS-short-curve}, there exists a constant $\epsilon>0$
depending only on $r_{\max}$ such that for any simple loop
$C\subset\partial\Nscr_{r}(\gamma)$ with $\Length(C)<2\pi\sinh{}r$,
each least area disk $\Delta\subset\Nscr_{r}(\gamma)$ bounded $C$ is
$\epsilon$ away from $\gamma$, i.e.
$\dist(\Delta,\gamma)\geq\epsilon$, where
$r\in[r_{\max}/2,r_{\max}]$. Next choose $t$ sufficiently large such
that $\sigma(1-t)<\epsilon$.

For any $t\in[0,1)$, $S_{t}\cap\T(\gamma)$ is either empty or
the union of annuli and  disks,
otherwise $S_{t}$ would be not incompressible.
For any sufficiently large $t\in[0,1)$, if
$S_{t}\cap\Nscr_{r_{\max}/2}(\gamma)=\emptyset$, then we are
done; so we assume
\begin{equation*}
   S_{t}\cap\Nscr_{r_{\max}/2}(\gamma)\ne\emptyset\ ,
   \quad\text{for all large}\ t\in[0,1)\ .
\end{equation*}
We don't need worry about the null-homotopic loops on
$\partial\Nscr_{r}(\gamma)$ with (hyperbolic) length $<2\pi\sinh{}r$
for all $r\in[r_{\max}/2,r_{\max}]$. In fact, suppose
$\Delta\subset{}S_{t}\cap\Nscr_{r}(\gamma)$ is a disk such that
$\Length(\partial{}\Delta)<2\pi\sinh{}r$, then $\partial{}\Delta$
bounds a least area minimal surface
$\Delta'\subset\Nscr_{r}(\gamma)\setminus\Nscr_{\epsilon}(\gamma)$
with respect to the hyperbolic metric, which is also a $g_{t}$-least
area disk since $\sigma(1-t)<\epsilon$, this implies $\Delta=\Delta'$.
So $\Delta$ is disjoint from $\Nscr_{\epsilon}(\gamma)$.

Now we can assume that for any $r\in[r_{\max}/2,r_{\max}]$,
the boundary of each component of
$S_{t}\cap\Nscr_{r}(\gamma)$ which is a disk has
length $>2\pi\sinh{}r$. By
Lemma \ref{lem:CG2006-(1.18)} and coarea formula
\eqref{eq:coarea formula II}, we have estimate
\begin{align*}
   \Area(S_{t}\cap\T(\gamma),g_{t})
      &\geq\Area(S_{t}\cap\T(\gamma))\\
      &\geq\Area(S_{t}\cap(\T(\gamma)\setminus
           \Nscr_{r_{\max}/2}(\gamma)))\\
      &\geq\int_{r_{\max}/2}^{r_{\max}}
           \Length(S_{t}\cap\partial\Nscr_{s}(\gamma))\,ds\ ,
\end{align*}
here $\Area(\cdot,g_{t})$ means $g_{t}$-area and
$\Area(\cdot)$ means hyperbolic area. There are two cases:
\begin{enumerate}
   \item if $S_{t}\cap\Nscr_{s}(\gamma)$ has a
         component that is an annulus, then
         \begin{align*}
            \Length(S_{t}\cap\partial\Nscr_{s}(\gamma))
               &\geq\Length(\Acal_{\theta/l}
                    \cap\partial\Nscr_{s}(\gamma))\\
               &=2(l^{2}\cosh^{2}s+
                 \theta^{2}\sinh^{2}s)^{1/2}\\
               &\geq{}2\theta\sinh{}s\ ,
         \end{align*}
   \item if $S_{t}\cap\Nscr_{s}(\gamma)$ has a
         component that is a disk, then by assumption,
         \begin{equation*}
            \Length(S_{t}\cap\partial\Nscr_{s}(\gamma))\geq
            2\pi\sinh{}s\ .
         \end{equation*}
\end{enumerate}
When $\theta/l\gg{}1$, so $r_{\max}\gg{}1$ and
$l\sinh{}r_{\max}\ll{}1$, then we have the following
estimates
\begin{equation*}
   \Area(S_{t}\cap\T(\gamma))
   \geq{}2c(\theta)\int_{r_{\max}/2}^{r_{\max}}\sinh{}s\,ds
   >{}c(\theta)\cosh{}r_{\max}
\end{equation*}
where $c(\theta)=\min\{\theta,\pi\}$, and
\begin{equation*}
   \Area(\partial\T(\gamma))
   =2\pi{}l\sinh{}r_{\max}\cosh{}r_{\max}
   <{}c(\theta)\cosh{}r_{\max}\ .
\end{equation*}
Recall that $g_{t}$ metric is equal to the hyperbolic metric
outside $\Nscr_{\sigma(1-t)}(\gamma)$ and
$\sigma(1-t)<\epsilon$, thus we have
\begin{equation*}
   \Area(S_{t}\cap\T(\gamma), g_{t})>
   \Area(\partial\T(\gamma),g_{t})\ .
\end{equation*}
Now we can use cutting and pasting, i.e. replace each
component of $S_{t}\cap\T(\gamma)$ by either an annulus or
a disk on $\partial\T(\gamma)$, to get a new surface
$S_{t}'\subset{}M\setminus\Nscr_{\epsilon}(\gamma)$ having 
the following properties
\begin{itemize}
   \item $\Area(S_{t}',g_{t})<\Area(S_{t},g_{t})$,
   \item $S_{t}'\cap\Nscr_{\epsilon}(\gamma)=\emptyset$, and   
   \item $S_{t}'$ is isotopic to $S$ in $M\setminus\gamma$.
\end{itemize}
So $S_{t}$ must be disjoint from $\Nscr_{\epsilon}(\gamma)$
for all sufficiently large $t$.
\end{proof}

By the Claim above, we see that $S_{t}$ is a least area minimal
surface with respect to the {\em hyperbolic metric} and
$S_{t}\cap\Nscr_{\epsilon}(\gamma)=\emptyset$, where $\epsilon$
is a constant independent of $t$ for all sufficiently large $t$,
therefore all of them must be
the same surface when $t$ is sufficiently large. As $t\to{}1$,
we get a least area minimal surface $T$ satisfying the conditions
list in this theorem.
\end{proof}

\begin{corollary}\label{cor:corollary of main theorem}
Let $\Gamma=\{\gamma_{i}\}_{i=1}^{n}\subset{}M$ be the 
collection of mutually disjoint simple closed geodesics, with 
complex length $\Lscr_{i}=l_{i}+\sqrt{-1}\,\theta_{i}$, such 
that $\theta_{i}/l_{i}\gg{}1$ for $1\leq{}i\leq{}n$,
and let $S\subset{}M\setminus\Gamma$ be a closed embedded 
incompressible surface in $M$, then there exists a least area 
minimal surface $T\subset{}M\setminus\Gamma$ such that $T$ 
is isotopic to $S$ in $M\setminus\Gamma$.

In particular, $M$ contains at least $2^{n}$ minimal surfaces,
which are incompressible, and which are not isotopic to each
other in $M\setminus\Gamma$. 
\end{corollary} 

\begin{proof}By Theorem \ref{thm:meyerhoff1987-(7)}, the tubes
$\T(\gamma_{1}),\T(\gamma_{2}),\ldots,\T(\gamma_{n})$ are 
mutually disjoint, so replace $\gamma$ by $\Gamma$ in the 
proof of Theorem \ref{thm:main theorem}, we can prove the 
first part of the corollary in the same way.

For the second part, J.-P. Otal proved in 
\cite[Theorem B]{otal2003} that $\Gamma$
is {\em unlinked} in the following sense: There exists a 
homeomorphism between $M$ and $\Sigma\times\R$ 
such that each component of $\Gamma$ is contained in one 
of the surfaces $\Sigma\times{i}$, $1\leq{}i\leq{}n$. Therefore
the (incompressible) surfaces in $M\setminus\Gamma$ can 
separate $\Gamma$ in
\begin{equation*}
   1+\binom{n}{1}+\binom{n}{2}+\cdots+\binom{n}{n-1}+1=
   2^{n}
\end{equation*}
ways such that they are not isotopic to each other in
$M\setminus\Gamma$. By the first part, we have at least
$2^{n}$ minimal surfaces.
\end{proof}

\bibliographystyle{amsalpha}
\bibliography{ms(QF)ref}

\end{document}